\documentclass[12pt]{amsart}
\usepackage{amssymb,latexsym}
\usepackage{enumerate}

\makeatletter
\@namedef{subjclassname@2010}{%
  \textup{2010} Mathematics Subject Classification}
\makeatother

\newtheorem{thm}{Theorem}[section]

\theoremstyle{definition}

\numberwithin{equation}{section}

\newtheorem{theorem}{Theorem}[section]

\theoremstyle{definition}
\newtheorem{definition}[theorem]{Definition}

\theoremstyle{remark}

\begin{document}


\baselineskip=17pt



\title[New examples of divisibility sequences]{New examples of determinant divisibility sequences}
\author[K. G\'ornisiewicz]{Krzysztof G\'ornisiewicz}
\address{Faculty of Mathematics and Computer Science\\ Adam Mickiewicz University\\
61-614 Pozna\'n, Poland}
\email{krisgorn@amu.edu.pl}
\thanks{Grant: UMO-2012/07/B/ST1/03541}
\subjclass[2010]{Primary 15A18, 14E20; Secondary 11A51, 11B39 , 11R32}
\keywords{divisibility sequence, eigenvalues, Jacobian matrix, algebraic numbers}

\begin{abstract}
In this paper  we consider divisibility sequences obtained from square matrices. We work with of matrix divisibility sequences associated to a semigroup  and  arising from endomorphisms of an affine space.  We prove that determinant divisibility sequences originated from  powers of square matrices are generalized Lucas sequences.
\end{abstract}

\maketitle

\section{Introduction}

By the divisibility sequence we mean in this paper  a sequence  $\{d_n\}_{n \in \mathbb{N}}$ of integers such that  if $ n|m $ then $d_n|d_m$. One of the most famous divisibility sequence is the Fibonacci  sequence: 0, 1,  1,  2, 3, 5, 8, 13, 21, 34,... which arise from linear recurrence:  $F_n=F_{n-1}+F_{n-2}$. This is an example of  the Lucas sequences: $L_n=\frac{\alpha^n-\beta^n}{\alpha-\beta}$, where $\alpha,\beta$ are the roots of some quadratic polynomial over $\mathbb{Z}$. See [2] for a complete classification of linear recurrence divisibility sequences and [5], [6] for introduction to other divisibility sequences.
In this paper we discus properties of certain matrix divisibility sequences. We follow the approach initiated in [1].\\

\section{Matrix divisibility sequence}
Let $S$ be a commutative ring with 1.  Let $M_r(S)$ be  a  ring of  $r \times r$ matrices with entries in $S$.  By a divisor class of a matrix $M \in M_r(S)$ we mean a coset $GL_r(S) \cdot  M$ of $M$ with respect to the natural left action of $GL_r(S)$ . We say that matrix $M\in M_r(S)$  divides a matrix $N\in M_r(S)$ if there exists a  $Q\in M_r(S)$ such that $N = QM$. If $M$ divides $N$, then any element of the divisor class of $M$ also divides $N$.  Let $(\Gamma, \cdot)$ denote a semigroup. A \textit{divisibility sequence of matrices} over a commutative ring S, indexed by $\Gamma$, is a collection of matrices $\{M_{\alpha}\}_{\alpha\in\Gamma}$ in $M_r(S)$, such that if $\alpha$  divides $\beta$ in $\Gamma$, then $M_{\alpha}$ divides $M_{\beta}$ in $M_r(S)$. 
If $\{M_{\alpha}\}_{\alpha\in\Gamma}$ is a divisibility sequence of matrices, then  by the  multiplicativity of the determinant $\{det(M_{\alpha})\}_{\alpha\in\Gamma}$ is a divisibility sequence of elements of the ring S.\\ 
We fix  a faithful representation:
$$
[\cdot]: \Gamma \hookrightarrow \mbox{End}(\mathbb{A}_S^m):\alpha \rightarrow [\alpha]
$$
of $\Gamma$ into the group of endomorphisms of affine m-dimensional space $\mathbb{A}_S^m$ over $S$.\\

\begin{definition}
Let $x \in \mathbb{A}^r_S$. The \textit{matrix divisibility sequence} associated to  $(\Gamma,[\cdot])$ is the sequence of Jacobians $\{J_{\alpha}(x)\}_{\alpha\in\Gamma}$ which are $r\times r$ matrices with $(i,j)-$entry given by partial differentials:
$$
[J_{\alpha}(x)]_{i,j}:=\partial(([\alpha](x))_i)/ \partial x_j,
$$
where $([\alpha](x))_i$ is an  ith entry of the value of the endomorphism $[\alpha]$ on $x$.
The associated \textit{determinant divisibility sequence} is defined by
$
\{det(J_{\alpha}(x)\}_{\alpha\in\Gamma}.
$\\
\end{definition}

\section{Main Result}

\begin{thm}
Let $X \in GL_r(\mathbb{Z})$ and $\lambda_1,\ldots,\lambda_r$ be eigenvalues of $X$. Then for every $n \geq 1$:
\begin{equation}
D_n=n^2\left[detX\right]^{n-1}
 \prod_{1\leq i< j\leq r}\left(\frac{\lambda_i^n-\lambda_j^n}{\lambda_i-\lambda_j}\right)^2  \quad\quad 
\label{Dn}
\end{equation}
is an integer and the sequence $\{D_n\}_{n\in\mathbb{N}}$ is a determinant divisibility sequence.\\
\end{thm}

\noindent
\textbf{Proof}: Let $X, Y, Z$ be square $s\times s$ matrices. Assume that entries of matrices $ Y$ and $Z$ are functions of entries of the matrix $X$. Then  the following matrix derivative formula holds ([4]):
\begin{equation}
\frac{d(YZ)}{dX}=(I\otimes Y)\frac{dZ}{dX}+(Z^t\otimes I)\frac{dY}{dX}, \quad\quad 
\label{md}
\end{equation}
where $\otimes$ means the Kronecker product, $I$ is the identity matrix of rank $s$ and $A^{t}$ means the transpose matrix of $A$. 
In addition we will use property of the Kronecker product:
\begin{equation}
(A\otimes C)(B\otimes D)=AB\otimes CD\,\,  \label{pK}
\end{equation}
for any square matrices $A, B, C, D$ of size $s\times s$. 
From now on we fix $\Gamma=\mathbb{N}$. Consider the group $G$ of all invertible $s\times s$ matrices with the embedding:\\
$$
G \rightarrow \mathbb{A}^{s^2}:
\left[\begin{array}{ccc}
X_{11}&\cdots&X_{1s}\\ \vdots&\ddots&\vdots\\X_{s1}&\cdots&X_{ss}
\end{array}\right]
\mapsto
\left(X_{11},\ldots,X_{1s},\ldots,X_{s1},\ldots,X_{ss}\right)
$$
We define the endomorphism $[n]$ for $n\in\mathbb{N}$. Let $X:= \left[X_{ij}\right] \in G$ and respectively $X^n:=[\bar{X}_{kl}]\in G$, where we treat $\bar{X}_{kl}$  as functions of $X_{ij}$, for $1\leq i,j,k,l\leq s$. We define $[n]:\mathbb{A}^{s^2}\rightarrow \mathbb{A}^{s^2}$ as $$[n]\left(X_{11},\ldots,X_{1s},\ldots,X_{s1},\ldots,X_{ss}\right)=
\left(\bar{X}_{11},\ldots,\bar{X}_{1s},\ldots,\bar{X}_{s1},\ldots,\bar{X}_{ss}\right).$$
Using (\ref{md}) we compute Jacobians of the $n$-th power of the matrix $X$\\
$$
J_n=\frac{d(X^n)}{dX}=\frac{d(X^{n-1}X)}{dX}=(I\otimes X^{n-1})\frac{dX}{dX}+(X^t\otimes I)\frac{d X^{n-1}}{dX}.$$
By induction and  the property (\ref{pK}) of the Kronecker product we get:
$$
J_n=\sum^{n-1}_{k=0}(X^t)^k\otimes X^{n-1-k}.
$$

\noindent
Assume that $X\in GL_s(\mathbb{Z})$ is a diagonalizable matrix  . Then $X=PDP^{-1}$, for some $P\in \mbox{GL}_s(\mathbb{C})$ and a diagonal matrix $D= diag\{\lambda_1,\lambda_2,\dots,\lambda_s\}$.  The set of matrices with distinct eigenvalues  is dense in the set of all square matrices $M_s(\mathbb{C})$ with respect to the topology of $\mathbb{C}^{s^2}$. Hence we can assume that $X$ has different eigenvalues.  Therefore
$$
J_n=\sum^{n-1}_{k=0}((PDP^{-1})^t)^k\otimes (PDP^{-1})^{n-1-k}=
\sum^{n-1}_{k=0}(P^{-1})^t (D^k)^t(P)^t\otimes PD^{n-1-k}P^{-1}=
$$
$$
=((P^{-1})^t\otimes P)\left[\sum^{n-1}_{k=0}(D^k\otimes D^{n-1-k})\right](P^t\otimes P^{-1})
$$
and the determinant divisibility sequence is of the form:
$$
D_n=\det J_n=\det\left[\sum^{n-1}_{k=0}(D^k\otimes D^{n-1-k})\right].
$$
Since $D^k\otimes D^{n-1-k}$ is a diagonal matrix whose diagonal consists of terms $\lambda_i^k\lambda_j^{n-1-k}$, we conclude that:
$$
D_n=n^2\prod_{l=0}^s\lambda_l^{n-1}\prod_{1\leq i\neq j\leq s} \sum_{k=0}^{n-1}\lambda_i^k\lambda_j^{n-1-k}=n^2\left[detX\right]^{n-1}
 \prod_{1\leq i< j\leq s}\left(\frac{\lambda_i^n-\lambda_j^n}{\lambda_i-\lambda_j}\right)^2.
$$ 
The last term on the right hand side is a product of values of symmetric polynomials computed at eigenvalues of the matrix $X$.  The Galois group of the splitting field of the characteristic polynomial of $X$ acts trivially on these algebraic integers, hence  $D_n \in \mathbb{Z}$.  For any  $n,m$ such that $n$ divides $m$ we have $[\det X]^n | [\det X]^m$ and $(\lambda_i^n-\lambda_j^n)|(\lambda_i^m-\lambda_j^m)$ .  Therefore, the sequence $D_n$ is a divisibility sequence. 
For equal eigenvalues we compute $D_n$  using the exact form of symmetric polynomials instead of their fractional expression.

\noindent
If $X$ is not  a diagonalizable matrix, then instead of the matrix $D$ we consider an upper-triangular matrix obtained from the Jordan form of $X$ with eigenvalues $\lambda_1,\lambda_2,\dots,\lambda_r$ on the diagonal  and eventually we get the same sequence $D_n$.

\section{Examples}

\bigskip
\noindent
1) Let $X\in \mbox{GL}_2(\mathbb{Z})$ and $a=\mbox{tr}X, b=\mbox{tr}^2X-4\mbox{det}X$. Then using Theorem 1 we obtain the sequence  presented in [1], example 4.3:
$$
D_n=\frac{n^2}{b}\left[\mbox{det}X\right]^{n-1}\left(\left( \frac{a+\sqrt{b}}{2}\right)^n- 
\left( \frac{a-\sqrt{b}}{2}\right)^n\right)^2
$$
\noindent
2) Let  $X\in \mbox{GL}_3(\mathbb{Z})$ and  $b =  -\mbox{tr}X, \:\:c=X_{11}+X_{22}+X_{33}, \:\:d=-\det X $. The discriminant of the characteristic polynomial of $X$  is $\Delta=(4\Delta_0^3-\Delta_1^2)/27$, where   $\Delta_0=b^2-3c$ and $\Delta_1=2b^3-9bc+27d$.  We obtain the divisibility sequence defined by:
$$
D_n=\frac{n^2d^{n-1}}{\Delta}\prod_{i=1}^3\left[\left(\frac{b+\epsilon^i A+\epsilon^{2i}\Delta_0\bar A}{3}\right)^n-\left(\frac{b+\epsilon^{i+1} A+\epsilon^{2i+2)}\Delta_0\bar A}{3}\right)^n\right]^2,
$$
where $A=\sqrt[3]{(\Delta_1+\sqrt{-27\Delta})/2}$, $\bar A=\sqrt[3]{(\Delta_1-\sqrt{-27\Delta})/2}$ and $\epsilon$ is a fixed primitive cube root of unity.\\
\noindent
3) It is easy to compute values of $D_n$ for any square matrix $X$. The matrix $X=\left[\begin{array}{rrr}
1&-2&-6\\0&1&3\\-1&0&1
\end{array}\right]\in Gl_3(\mathbb{Z})$ 
gives the following divisibility sequence: \\

\bigskip

\setlength{\tabcolsep}{10pt}.
\renewcommand{\arraystretch}{1.3}
\begin{tabular}{|l|l|l|} 
  	 \hline
$n$&$d_n$&factorization of $d_n$\\
	\hline
    $1$ & $1$ &1\\
    $2$&$100$&$2^{2}5^{2}$\\
    $3$&$6561$&$3^{8}$\\
    $4$&$193600$&$2^{6}5^{2}11^{2}$\\
    $5$&$808201$&$29^{2}31^{2}$\\
    $6$&$189612900$&$2^{2}3^{8}5^{2}17^{2}$\\
    $7$&$50131657801$&$41^{2}43^{2}127^{2}$\\
    $8$&$4096576000000$&$2^{12}5^{6}11^{2}23^{2}$\\
    $9$&$159625511221401$&$3^{14}53^{2}109^{2}$\\
    $10$&$1865976489302500$&$2^{2}5^{4}29^{2}31^{6}$\\
    $11$&$31583922467632921$&$131^{2}857^{2}1583^{2}$\\
    $12$&$21985833099924302400$&$2^{6}3^{8}5^{2}11^{2}17^{2}71^{2}109^{2}$\\
    $13$&$2370466451421685365841$&$1637^{2}4057^{2}7331^{2}$\\
    $14$&$118070682478980566428900$&$2^{2}5^{2}41^{2}43^{6}83^{2}127^{2}$\\
    $15$&$2362255369723766871090801$&$3^{8}29^{2}31^{2}2969^{2}7109^{2}$\\
    $16$&$84956038709284864000000$&$2^{18}5^{6}11^{2}23^{2}47^{2}383^{2}$\\
	\hline
\end{tabular}\\

\noindent
4) The matrix $X=\left[\begin{array}{rrrr}
-1&2&4&-1\\0&1&-2&2\\-1&0&-1&0\\0&1&0&1
\end{array}\right]\in Gl_4(\mathbb{Z})$ 
gives the following divisibility sequence: \\

\setlength{\tabcolsep}{10pt}.
\renewcommand{\arraystretch}{1.3}
\begin{tabular}{|l|p{5.2cm}|l|} 
  	 \hline
$n$&$d_n$&factorization of $d_n$\\
	\hline
    $1$&$1$ &1\\ \hline
    $2$&$65536$&$2^{16}$\\ \hline
    $3$&$1$&1\\ \hline
    $4$&$281474976710656$&$ 2^{48}$\\ \hline
    $5$&$18448995933652254721$&$4295229439^{2}$\\ \hline
    $6$&$18013780039499776$&$2^{16}7^{2}74897^{2}$\\ \hline
    $7$&$7922332684788105606123945$ $9841$&$281466386710529^{2}$\\ \hline
    $8$&$5194832314440011219064571$ $543158784$&$2^{60}7^{4}23^{2}59561^{2}$\\ \hline
    $9$&$5775028020578783636854257$ $0774529$&$37^{2}701^{2}292993041329^{2}$\\ \hline
    $10$&$2229666183092939997026658$ $7262959037621272576$&$ 2^{16}19^{4}3449^{4}4295229439^{2}$\\ \hline
    $11$&$1463330673647120201450844$ $900178197550156472647681$&$ 32363^{2}7282397^{2}5132726390881^{2}$\\ \hline
    $12$&$9132817257932632786870155$ $6304335790376407269376$&$ 2^{48}7^{2}13^{2}10177^{2}74897^{2}259691^{2}$\\ \hline
    $13$&$6274228310768040852924579$ $1977173633010223354344560$ $089620481$&
$3^{18}3769^{2}15053^{2}27205307^{2}2607270173^{2}$\\ \hline
    $14$&$4124060734576866740544330$ $9242707472643430878237415$ $8659336863744$&$
 2^{16}13^{2}794009^{2}27304061^{2}281466386710529^{2}$\\ \hline
    $15$&$9806175554643239147044270$ $0484791252942607177394577$ $588251525121$&$ 17489^{2}4295229439^{2}131825214490835791^{2}$\\ \hline
    $16$&$1765121615339370515604475$ $3103664126594001046686372$ $ 42611924881667927310336$&$
2^{72}7^{8}23^{2}59561^{2}20394769^{2}288208447^{2}$\\	\hline
\end{tabular}

\newpage
\textbf{References:}\\

[1] Gunther Cornelissen, Jonathan Reynolds, Matrix divisibility sequences, Acta Arith. 156 (2012), 177-188

[2] J.P. Bezivin, A. Ptheo, A.J. van der Poorten, A full characterization of divisibility sequences,
Amer. J. of Math.112 (1990), 985–1001;

[3] Rachel Shipsey, Elliptic divisibility sequences, Ph.D. thesis, Goldsmiths College, University of London, see home-pages.gold.ac.uk/rachel, 2000

[4] Mike Brooks: The Matrix Reference Manual,\\ http://www.ee.ic.ac.uk/hp/staff/dmb/matrix/intro.html

[5] Joseph H. Silverman,Generalized greatest common divisors, divisibility sequences, and Vojta's conjecture for blowups, Monatsh. Math. 145 (2005), no. 4, 333–350.

[6] Patrick Ingram, Elliptic divisibility sequences over certain curves, J. Number Theory 123 (2007), no. 2, 473–486.

\end{document}